\documentclass[journal]{IEEEtran}
\IEEEoverridecommandlockouts

\usepackage{amsmath}
\usepackage{amssymb}
\usepackage{amsthm}
\usepackage{multirow}		
\usepackage{subfig}
\usepackage{stfloats}
\usepackage{xfrac}			
\usepackage{booktabs}
\usepackage{color}

\usepackage{accents}
\newcommand{\ubar}[1]{\underaccent{\bar}{#1}}

\allowdisplaybreaks

\newtheorem{Thm}{Theorem}
\newtheorem{Def}{Definition}

\newtheorem{remark}{Remark}
\newcommand{\beq}{\begin{equation}}
\newcommand{\eeq}{\end{equation}}
\newcounter{l1}
\newcounter{l2}
\newcounter{l3}
\newcommand{\bdotlist}{\begin{list}{$\bullet$}{}}
\newcommand{\bboxlist}{\begin{list}{$\Box$}{}}
\newcommand{\bbboxlist}{\begin{list}{\raisebox{.005in}{{\tiny $\blacksquare$ \\ }}}{}}
\newcommand{\bdashlist}{\begin{list}{$-$}{} }
\newcommand{\blist}{\begin{list}{}{} }
\newcommand{\barablist}{\begin{list}{\arabic{l1}}{\usecounter{l1}}}
\newcommand{\balphlist}{\begin{list}{(\alph{l2})}{\usecounter{l2}}}
\newcommand{\bAlphlist}{\begin{list}{\Alph{l2}.}{\usecounter{l2}}}
\newcommand{\bdiamlist}{\begin{list}{$\diamond$}{}}
\newcommand{\bromalist}{\begin{list}{(\roman{l3})}{\usecounter{l3}}}

\begin{document}

\title{\huge{A General Economic Dispatch Problem with Marginal Losses}}

\author{M. Garcia\IEEEauthorrefmark{3}, R. Baldick\IEEEauthorrefmark{3}, S. Siddiqi\IEEEauthorrefmark{4} \vspace{-17pt}
\thanks{\IEEEauthorrefmark{3} Electrical and Computer Engineering, U.T. Austin}%
\thanks{\IEEEauthorrefmark{4} Crescent Power, Inc.}%
\thanks{This work was supported by the Defense Threat Reduction Agency.}
\thanks{Some of this work was accomplished during a visit in Fall 2017 by the second author to the Florence School of Regulation, European University Institute, Fiesole, Italy.  Support from the University of Texas Faculty Development Program and Professor Jean-Michel Glachant of the Florence School of Regulation is gratefully acknowledged.}

}

\maketitle


\begin{abstract}
Standard economic dispatch problems that consider line losses are linear approximations of a non-convex economic dispatch problem formulated by fixing voltage magnitudes and assuming the decoupling of real and reactive power.  This paper formulates and analyzes the general non-convex economic dispatch problem, incorporating and generalizing the Fictitious Nodal Demand (FND) model, resulting in a slack bus independent formulation that provides insight into standard formulations by pointing out commonly used but unnecessary assumptions and by deriving proper choices of ``tuning parameters.'' The proper choice of loss allocation is derived to assign half of the losses of each transmission line to adjacent buses, justifying approaches in the literature.  Line constraints are proposed in the form of voltage angle difference limits and are proven equivalent to various other line limits including current magnitude limits and mid-line power flow limits.  The formulated general economic dispatch problem with marginal losses consistently models flows and loss approximation, results in approximately correct outcomes and is proven to be reference bus independent. Various approximations of this problem are compared using realistically large transmission network test cases.\vspace{-10pt}

\end{abstract}\normalfont

\vspace{-2pt}
\section{Introduction}
\vspace{-3pt}
Independent System Operators (ISOs) governed by the Federal Energy Regulatory Commission (FERC) claim to realize significant benefits by incorporating marginal losses into economic dispatch \cite{hogan2017priorities}.  In fact, PJM has reported 100 million dollars of savings per year in energy and congestion costs~\cite{PJMTraining}.  Long term benefits may be realized through the accurate formation of Locational Marginal Prices (LMPs), which better guide investment and operational decisions.  These claims have encouraged other ISOs to consider implementing marginal losses into their economic dispatch~\cite{ICF2017}.  

Despite the increased prevalence of marginal losses, a standard economic dispatch problem has not emerged in industry, so that energy markets that consider marginal losses have varying implementations~\cite{eldridge2017marginal}.  The choice of implementation of marginal losses have both market price and resource dispatch impacts that may significantly benefit some players while significantly disadvantaging other players as compared to \emph{Alternating Current Optimal Power Flow} (AC OPF) outcomes.  To better understand these various implementations, this paper derives an economic dispatch problem from first principles that serves as a generalization of those used in practice, including the Fictitious Nodal Demand (FND) model from reference \cite{li2007dcopf}, and explains the various assumptions and approximations required to attain different practical formulations.  Common and unnecessary assumptions are identified and proper choices of tuning parameters are specified.  Certain approximations are shown to increase the error in price and dispatch outcomes.



The real-power loss of a transmission line can be roughly approximated as the product of the line resistance and the squared real-power flowing through the line~\cite{hogan2002financial,philpott1999experiments}.  Though this is not the most accurate approximation, it is intuitive and provides insight into why power balance equations that consider marginal losses are inherently non-linear and cause economic dispatch problems to be non-convex.  This non-convexity makes the problem difficult to solve in principle, and as a result ISOs typically approximate the non-linear loss term by linearizing around some base-case~\cite{coffrin2012approximating,yang2017lmp,yang2018linearized}.  Many works analyze these linear optimization problems without initially formulating the original non-convex problem~\cite{eldridge2017marginal,litvinov2004marginal,Zhang2013,akbari2014linearized}.  These linearization techniques intend to approximate the solution to the original non-convex economic dispatch problem.  However, the combination of an initial linearization followed by an ad hoc incorporation of an approximation to nonlinear losses makes it unclear how to evaluate the approximation errors in these formulations.  Moreover, such an approach typically produces a dependence of approximations and therefore approximation errors on the choice of Loss Distribution Factors (LDFs), which we wish to avoid. To better understand the approximation errors, and avoid dependence on LDFs, this paper focuses on formulating the fundamental non-convex economic dispatch problem and several nonlinear approximations, leaving linearization techniques to be addressed in future work.  In comparison to its linearized counterpart, the nonlinear approximations are more difficult to solve. However, we illustrate that computational efforts are modest and global optimality can often be guaranteed on realistically large test cases providing a more accurate approximation to the exact AC OPF problem.

Economic dispatch models that consider marginal losses typically use variations of the DC model that represent losses as load allocated to one or more buses throughout the system.  See reference \cite{stott2009dc} for an in depth description of variations of the DC model of the transmission system.  Reference \cite{stott2009dc} restricts its analysis to DC models with constant nodal loss allocation.  However, accurate price setting requires losses to be represented, not as constant, but as a varying function of the system state.  With this in mind, economic dispatch formulations should accurately capture the approximately quadratic relationship between losses and real-power flow.  

Many previous works augment a lossless DC power flow model with a quadratic loss model, eg. \cite{li2007dcopf,li2011fully}. As suggested above, this approach is somewhat self-contradictory because the DC model of the system is only accurate when lines are lossless. Our paper will shed light on why such an approach nevertheless results in a reasonable approximation to overall losses. However, by assuming lossless transmission lines it is not clear where the additional load due to losses should be allocated.  Initial formulations allocated all losses to the slack bus \cite{shahidehpour2002market}, which results in at least some of the losses being modeled as occurring far from the lines that actually incur the losses.  Reference \cite{litvinov2004marginal} recognized that the solution to the resulting dispatch problem depended on the choice of slack bus and corrected this problem by introducing ``Loss Distribution Factors'' that fix the fraction of total system losses allocated to each bus.  However, the solution to the dispatch problem in turn depends on the choice of those LDFs and the authors do not provide a method for determining these factors and do not clarify how to remove this fundamental dependence on the choice of LDFs.  More accurate models include the Fictitious Nodal Demand (FND) model from references \cite{li2007dcopf} and \cite{hogan2002financial}, which allocate half of the losses of each line to both adjacent buses.  Our work, which builds on the FND model, provides insight into how losses should be allocated as fictitious demand and does not require designation of a slack bus or LDFs.  We acknowledge that without use of LDFs our formulation must explicitly represent the loss allocation as a variable for each bus in the system, as opposed to representing the total system losses as a single variable, increasing the dimension of the feasible set and making the problem more difficult to solve.

Our work formulates an economic dispatch problem that generalizes the FND problem formulation by removing the standard DC assumptions in its derivation and using only the key assumption that voltage magnitudes are fixed. That is, for the first time, we rigorously show that the standard FND derivation of loss approximations results in an approximately correct answer despite the logical inconsistency of first ignoring losses and then re-incorporating losses as an ad hoc adjustment.  By formulating the flow and loss approximation consistently, our formulation is independent of the choice of reference and slack bus and provides insight into the relationship between the economic dispatch problem solved in practice and the AC OPF problem that includes real and reactive power as decision variables.  Specifically, dispatch problems used in practice serve as a linear approximation to the non-convex dispatch problem formulated in this paper, which in turn serves as an approximation to the AC OPF problem. This relationship will be outlined without directly formulating the AC OPF problem, which is detailed in other work \cite{cain2012history}. We emphasize that none of the formulations considered in this paper consider contingency constraints, which is an issue that should be addressed in future work.

The rest of the paper is organized as follows. Section \ref{sec:LineModel} provides a model of a transmission line using the assumption that voltage magnitudes are fixed.  We explain how this model is a generalization of other models used in the literature and outline various common approximations.  Section \ref{sec:OPFForm} uses a general transmission line model to formulate a non-convex economic dispatch problem that utilizes line constraints in the form of voltage angle difference limits across each transmission line.  This section continues to explain how the voltage angle limit parameters can be chosen to indirectly enforce limits on various physical parameters associated with a transmission line. Section \ref{sec:Results} provides empirical results analyzing the error of the various approximations as well as their relation to the AC OPF problem. Section \ref{sec:Conc} concludes.
\vspace{-2pt}

\vspace{-5pt}
\section{Transmission Line Model} \label{sec:LineModel}
\vspace{-3pt}
We begin by introducing notation.  Lower case subscripts are used to indicate elements of matrices/vectors.  For example the element in the $i^{th}$ row and $j^{th}$ column of matrix $M$ is denoted $M_{ij}$. The set of $n$ dimensional real numbers is denoted $\mathbb{R}^n$.  Furthermore, the system is modeled as an undirected graph $\mathcal{G}=(\mathbb{N},\mathbb{E})$ where $\mathbb{N}$ is the set of nodes (buses) and $\mathbb{E}$ is the set of edges (transmission lines). There are $n$ buses and $m$ transmission lines.  Lines will be indexed by $\ell\in \{1,\ldots,m\}$.  The imaginary number is denoted $\mathbf{i}$.  The transpose of a vector or matrix is indicated by a superscript dagger, eg. $M^{\dag}$.

\begin{figure}[h!]
\vspace{-10pt}
\begin{center}
\includegraphics[scale=1]{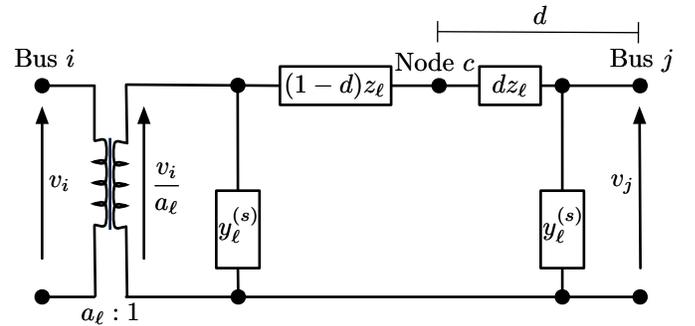}
\end{center}
\vspace{-5pt}
\caption{\label{PiModel} Circuit diagram of an arbitrary line $\ell$ connecting bus $i$ to bus $j$.}
\end{figure}

Figure \ref{PiModel} provides a circuit diagram of the general $\Pi$-model of a transmission line.  The transmission line is indexed by $\ell$ and connects bus $i\in \mathbb{N}$ to bus $j\in \mathbb{N}$.  The series impedance of the line is divided into two parts and separated by an intermediate node $c$ located at a fractional distance $d$ from from bus $j$. The total series impedance is denoted $z_{\ell}=r_{\ell}+\mathbf{i} x_{\ell}$ where $r_{\ell}$ is the series resistance and $x_{\ell}$ is the series reactance. The series impedance separating node $c$ from bus $j$ is in the amount $dz_{\ell}$.  We assume that the shunt conductances are zero, ie. the real part of $y^{(s)}_{\ell}$ is zero, and thus no real power flows through the shunt elements.  The complex voltage at bus $i$ is denoted $v_i=V_i \angle (\theta_i)$, where $V_i$ is the voltage magnitude and $\theta_i$ is the voltage angle.  An ideal transformer is located near bus $i$ with complex off-nominal turns ratio of $a_{\ell}= \tau_{\ell} \angle (\psi_{\ell})$.

\vspace{-5pt}
\subsection{Real-Power Flow on a Transmission Line}\label{TransMod}
\vspace{-2pt}
We consider the flow of real-power in the series element of the equivalent $\Pi$-model of a transmission line.  Assuming fixed voltage magnitudes, the power flowing through node $c$ in the series element in the direction of bus $j$ at an arbitrary fractional distance $d$ is represented by the following function:
\begin{align}
	\hat{F}_{\ell}(\Delta\!\theta_{\ell},d):=&  g_{\ell}({\textstyle\frac{d}{\tau_{\ell}^2}}V_i^2-(1-d)V_j^2)-{\textstyle\frac{b_{\ell}}{\tau_{\ell}}}V_iV_j\sin(\Delta\!\theta_{\ell}-\psi_{\ell}) \nonumber\\
	&\hspace{15pt}-{\textstyle\frac{g_{\ell}}{\tau_{\ell}}}V_iV_j(2d-1)\cos(\Delta\!\theta_{\ell}-\psi_{\ell}),\label{GeneralPowerFlow}
\end{align}
where $y_{\ell}=1/z_{\ell}=g_{\ell}+\mathbf{i}b_{\ell}$ and $\Delta\!\theta_{\ell}=\theta_i-\theta_j$. This model is a generalization of that in reference \cite{stott2009dc}. Notice that this function requires knowledge of the fixed voltage magnitudes.  Reference \cite{stott2009dc} suggests different ways of choosing the fixed voltage magnitudes including using the state estimated values or a local minimizer of the AC OPF problem. 

\vspace{-2pt}
\subsection{Loss Function}
\vspace{-2pt}
From equation (\ref{GeneralPowerFlow}), the power flowing into the line from bus $i$ and from bus $j$ are respectively expressed as $\hat{F}_{\ell}(\Delta\!\theta_{\ell},1)$ and $-\hat{F}_{\ell}(\Delta\!\theta_{\ell},0)$. The loss function represents the real power loss across the line and is derived by summing these two values:
\begin{equation}\hat{L}_{\ell}(\Delta\!\theta_{\ell}):=g_{\ell}(V_j^2+{\textstyle\frac{1}{\tau_{\ell}^2}}V_i^2)-2{\textstyle\frac{g_{\ell}}{\tau_{\ell}}}V_iV_j\cos(\Delta\!\theta_{\ell}-\psi_{\ell}).\label{LossFuncDef}\end{equation}
\subsubsection{Approximating the Loss Function}\label{LossFunctionApproximations}
ISOs may desire a simpler quadratic approximation in order to utilize quadratic programming software. A very accurate approximation uses a third order Taylor expansion of the cosine function in equation (\ref{LossFuncDef}) around $\Delta\!\theta_{\ell}=\psi_{\ell}$.  This approximation results in a quadratic model of losses since the coefficient of the third order term is zero, with quartic error on the order of $\frac{g_{\ell}V_iV_j}{24\tau_{\ell}}\Delta\!\theta_{\ell}^4$:
\begin{equation}\hat{L}_{\ell}(\Delta\!\theta_{\ell})\approx g_{\ell}(V_j^2+{\textstyle\frac{1}{\tau_{\ell}^2}}V_i^2)-{\textstyle\frac{g_{\ell}}{\tau_{\ell}}}V_iV_j(2-(\Delta\!\theta_{\ell}-\psi_{\ell})^2).\label{LossFuncFirstApprox}\end{equation}
Note that this approximation does not rely on any assumption that the line resistance is small compared to the reactance.  Though very accurate, using this approximation in practice still requires some knowledge of how to fix the voltage magnitudes and tap ratios.  An ISO may desire a simpler model that fixes voltage magnitudes to 1 p.u. and tap ratios to $a_{\ell}=1$.  Fixing these values results in the following approximation:
\begin{equation}\hat{L}_{\ell}(\Delta\!\theta_{\ell})\approx g_{\ell}\Delta\!\theta_{\ell}^2.\label{approx1}\end{equation}
Perhaps the most commonly used approximation equates the line losses to the product of resistance of the line and squared DC power flow across the line. We can arrive at this approximation from equation (\ref{approx1}) by additionally assuming $r_{\ell}\ll x_{\ell}$: 
\begin{equation}\hat{L}_{\ell}(\Delta\!\theta_{\ell})\approx r_{\ell} ({\textstyle\frac{1}{x_{\ell}}}\Delta\!\theta_{\ell})^2.\label{AnotherApprox}\end{equation}
This final approximation of the loss function is by far the simplest and carries the interpretation that losses are equal to the resistance times the squared DC power flow where the DC power flow is given by $\frac{1}{x_{\ell}}\Delta\!\theta_{\ell}$.  This interpretation is often used naively without an understanding of the several approximations used to get to this point.

\vspace{-5pt}

\subsection{Fictitious Nodal Demand Representation}\label{Sec:FND}

\begin{figure}[h!]
\vspace{-12pt}
\begin{center}
\includegraphics[scale=1]{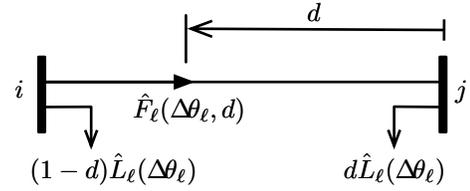}
\end{center}
\vspace{-15pt}
\caption{\label{FNDFig} One-line diagram of the FND representation of the proposed transmission line model.}
\vspace{-0pt}
\end{figure}

An FND representation of our transmission line model can be derived from equations (\ref{GeneralPowerFlow}) and (\ref{LossFuncDef}).  Notice that we can express the power flowing into the line from adjacent buses as follows, for an arbitrary fractional distance $d$:
\begin{align}
\hat{F}_{\ell}(\Delta\!\theta_{\ell}, 1) =&  \hat{F}_{\ell}(\Delta\!\theta_{\ell}, d) + (1-d)\hat{L}_{\ell}(\Delta\!\theta_{\ell}),\label{FNDEQ1}\\
-\hat{F}_{\ell}(\Delta\!\theta_{\ell}, 0) =&  -\hat{F}_{\ell}(\Delta\!\theta_{\ell}, d) + d \hat{L}_{\ell}(\Delta\!\theta_{\ell}).\label{FNDEQ2}	
\end{align}
These expressions lead to an FND model depicted in figure \ref{FNDFig} that is similar to that shown in reference \cite{stott2009dc}.  Specifically, a lossless transfer of power from bus $i$ to bus $j$ occurs in the amount $\hat{F}_{\ell}(\Delta\!\theta_{\ell}, d)$.  Losses are then represented as fictitious demand at bus $i$ in the amount $(1-d)\hat{L}_{\ell}(\Delta\!\theta_{\ell})$ and at bus $j$ in the amount $d\hat{L}_{\ell}(\Delta\!\theta_{\ell})$.  This FND model is equivalent to (\ref{GeneralPowerFlow}) in the sense that the buses see the same net power injection.

\subsection{Mid-line Power Flows}
A crucial observation is that $d$ can be chosen to simplify the form of the resulting model.  When $d=\frac{1}{2}$ the cosine term in equation (\ref{GeneralPowerFlow}) drops out, simplifying the expression for lossless power flow in the FND model (\ref{FNDEQ1}) and (\ref{FNDEQ2}).  This value is termed the \emph{mid-line power flow} and is written as follows:
\begin{equation}\hat{F}_{\ell}\left(\Delta\!\theta_{\ell},{\textstyle\frac{1}{2}}\right)={\textstyle\frac{1}{2}}g_{\ell}({\textstyle\frac{1}{\tau_{\ell}^2}}V_i^2-V_j^2)-{\textstyle\frac{b_{\ell}}{\tau_{\ell}}}V_iV_j\sin(\Delta\!\theta_{\ell}-\psi_{\ell}).\label{midlinefloweq}\end{equation}
We will see that with this choice of $d$, a linear approximation to $\hat{F}_{\ell}$ is accurate to second order in $\Delta\!\theta_{\ell}$.  This choice of $d$ is also convenient because of the symmetry of losses that occur around the mid-point of the line.  Specifically, the second term on the right-hand side of equations (\ref{FNDEQ1}) and (\ref{FNDEQ2}) are identically $\frac{1}{2}\hat{L}_{\ell}(\Delta\!\theta_{\ell})$.  This expression for the mid-line power flow is also in reference \cite{stott2009dc}.

\subsubsection{Approximating Mid-line Power Flows}\label{MidLinePowerFlowAproximations}
Unfortunately, the exact expression for the mid-line power flow (\ref{midlinefloweq}) is not convex in the vicinity of $\Delta\!\theta_{\ell}\approx 0$. However, a very accurate linear approximation can be attained through a second order Taylor expansion at $\Delta\!\theta_{\ell}=\psi_{\ell}$.  This approximation results in a linear model since the coefficient of the second order term is zero, with cubic error on the order of $\frac{b_{\ell}V_iV_j}{6\tau_{\ell}}\Delta\!\theta_{\ell}^3$:
\begin{equation}\hat{F}_{\ell}\left(\Delta\!\theta_{\ell},{\textstyle\frac{1}{2}}\right)\approx {\textstyle\frac{1}{2}}g_{\ell}({\textstyle\frac{1}{\tau_{\ell}^2}}V_i^2-V_j^2)-{\textstyle\frac{b_{\ell}}{\tau_{\ell}}}V_iV_j(\Delta\!\theta_{\ell}-\psi_{\ell}).\label{midlineflowfirstapprox}\end{equation}
Fixing the voltage magnitudes to 1p.u. and off-nominal tap ratios to $a_{\ell}\hspace{-2pt}=\hspace{-2pt}1$ further simplifies the previous \mbox{approximation to}:
\begin{equation}\hat{F}_{\ell}\left(\Delta\!\theta_{\ell},{\textstyle\frac{1}{2}}\right) \approx -b_{\ell}\Delta\!\theta_{\ell}.\label{midlinefloweqtaylor}\end{equation}
A final approximation often utilized in practice assumes \mbox{$r_{\ell}\ll x_{\ell}$}, leading to the following expression:
\begin{equation}\hat{F}_{\ell}\left(\Delta\!\theta_{\ell},{\textstyle\frac{1}{2}}\right) \approx {\textstyle\frac{1}{x_{\ell}}}\Delta\!\theta_{\ell}.\label{midlinefloweqfinal}\end{equation}
With respect to the FND representation of our model, this approximation can be interpreted as the addition of lossless DC power flow and line losses distributed equally to adjacent buses.  In fact, this interpretation is also consistent with with the FND formulation from references \cite{li2007dcopf} and \cite{hogan2002financial}, which suggest half losses of each line be allocated to both adjacent buses.  Our derivation shows that the allocation resulting from the choice $d=\frac{1}{2}$ is particularly advantageous in that it results in a simple expression for the flow that is well approximated by a linear function.

\vspace{-5pt}
\subsection{Squared Current Magnitude}\label{sec:curmag}
The squared magnitude of the current flowing into the transmission line from bus $i$ and bus $j$ can be represented as functions of the voltage angle difference $\Delta\!\theta_{\ell}$ and will be denoted $\hat{\mathcal{I}}_{ij}(\Delta\!\theta_{\ell})$ and $\hat{\mathcal{I}}_{ji}(\Delta\!\theta_{\ell})$ respectively.  These functions are respectively written as follows:
\begin{align}
\hspace{-5pt}\hat{\mathcal{I}}_{ij}(\hspace{-1pt}\Delta\!\theta_{\ell}\hspace{-1pt})&\hspace{-2pt}:=\hspace{-2pt}{\textstyle\frac{|y_{\ell}|^2}{\tau_{\ell}^2}}\hspace{-2pt}\left({\textstyle\frac{\alpha_{\ell}^2}{\tau_{\ell}^2}} V_i^2\hspace{-2pt}+\hspace{-2pt}V_j^2\hspace{-3pt}-\hspace{-2pt}{\textstyle\frac{2\alpha_{\ell}}{\tau_{\ell}}} V_iV_j\cos(\hspace{-2pt}\phi_{\ell}\hspace{-2pt}-\hspace{-2pt}\psi_{\ell}\hspace{-2pt}+\hspace{-2pt}\Delta\!\theta_{\ell}\hspace{-2pt})\hspace{-2pt}\right)\hspace{-3pt},\hspace{-5pt}\label{ExactCurrent}
\end{align}
\vspace{-10pt}
\begin{align}
\hspace{-5pt}\hat{\mathcal{I}}_{ji}\hspace{-1pt}(\hspace{-1pt}\Delta\!\theta_{\ell}\hspace{-1pt})&\hspace{-3pt}:=\hspace{-2pt}{\textstyle|y_{\ell}|^2}\hspace{-2pt}\left(\hspace{-2pt}{\textstyle\frac{1}{\tau_{\ell}^2}} V_i^2\hspace{-3pt}+\hspace{-2pt}\alpha_{\ell}^2 V_j^2\hspace{-3pt}-\hspace{-2pt}{\textstyle\frac{2\alpha_{\ell}}{\tau_{\ell}}} V_iV_j\hspace{-1pt}\cos(\hspace{-2pt}\phi_{\ell}\hspace{-2pt}+\hspace{-2pt}\psi_{\ell}\hspace{-2pt}-\hspace{-2pt}\Delta\!\theta_{\ell}\hspace{-2pt})\hspace{-3pt}\right)\hspace{-3pt},\hspace{-6pt}\end{align}

where $\alpha_{\ell}$ and $\phi_{\ell}$ are defined to be the magnitude and angle of the complex number $z_{\ell}(y_{\ell}^{(s)}+y_{\ell})$ respectively, so that $z_{\ell}(y_{\ell}^{(s)}+y_{\ell})=\alpha_{\ell} \angle (\phi_{\ell})$.  Notice that this complex number is approximately $1$ because the shunt admittance is typically much smaller than the series admittance.  As a result $\alpha_{\ell}\approx 1$ and $\phi_{\ell}\approx 0$.   In the following, we consider $\hat{\mathcal{I}}_{ij}(\Delta\!\theta_{\ell})$. The function $\hat{\mathcal{I}}_{ji}(\Delta\!\theta_{\ell})$ can be handled similarly.

\subsubsection{Approximating the Squared Current Magnitude}
Similar approximations to those in sections \ref{LossFunctionApproximations} and \ref{MidLinePowerFlowAproximations} can be used for the squared current magnitude function. The first approximation uses a third order Taylor expansion of the cosine function in equation (\ref{ExactCurrent}) around $\Delta\!\theta_{\ell}=-\phi_{\ell}+\psi_{\ell}$.  This approximation results in a quadratic function, with quartic error on the order of $\frac{2\alpha_{\ell} |y_{\ell}|^2 V_iV_j}{24\tau^3}(\Delta\!\theta_{\ell}+\phi_{\ell}-\psi_{\ell})^4$:
\begin{align}
\hspace{-5pt}\hat{\mathcal{I}}_{ij}\hspace{-1pt}(\Delta\!\theta_{\ell})&\hspace{-2pt}\approx\hspace{-4pt} \left(\hspace{-2pt}{\textstyle\frac{\alpha_{\ell}^2}{\tau_{\ell}^2}} V_i^2\hspace{-3pt}+\hspace{-2pt}V_j^2\hspace{-3pt}-\hspace{-2pt}{\textstyle\frac{\alpha_{\ell}}{\tau_{\ell}}} V_iV_j\hspace{-2pt}\left(2\hspace{-2pt}-\hspace{-2pt}(\Delta\!\theta_{\ell}\hspace{-2pt}+\hspace{-2pt}\phi_{\ell}\hspace{-2pt}-\hspace{-2pt}\psi_{\ell})^2\right)\hspace{-3pt}\right)\hspace{-4pt}{\textstyle\frac{|y_{\ell}|^2}{\tau_{\ell}^2}}\hspace{-1pt} .\hspace{-5pt}\label{ISFirstApprox}
\end{align}
In the case where shunt admittance $y_{\ell}^{(s)}$ is negligible we have $\alpha_{\ell}\hspace{-1pt}=\hspace{-1pt}1$ and $\phi_{\ell}\hspace{-1pt}=\hspace{-1pt}0$.  Additionally fixing the voltage magnitudes to 1p.u. and the off-nominal tap ratio to $a\hspace{-1pt}=\hspace{-1pt}1$ simplifies the previous approximation as follows:
\begin{align}
\hspace{-5pt}\hat{\mathcal{I}}_{ij}(\Delta\!\theta_{\ell})&\approx|y_{\ell}|^2 \Delta\!\theta_{\ell}^2.\label{IStaylor}
\end{align}
In the case where $r_{\ell}\ll x_{\ell}$ we have $y_{\ell}\approx \mathbf{i}b_{\ell}$.  The approximation (\ref{IStaylor}) then simplifies to the following:
\begin{align}
\hspace{-5pt}\hat{\mathcal{I}}_{ij}(\Delta\!\theta_{\ell})&\approx ({\textstyle\frac{1}{x_{\ell}}}\Delta\!\theta_{\ell})^2.\label{ISfinal}
\end{align}
This final approximation of the squared current magnitude is the simplest and carries the interpretation that current magnitude is equal to magnitude of the DC power flow.  

\subsection{General Transmission Line Model}
We have now introduced multiple approximations of the loss function $\hat{L}_{\ell}(\cdot)$, the mid-line power flow function $\hat{F}_{\ell}(\cdot,\frac{1}{2})$, and the squared current magnitude function $\hat{\mathcal{I}}_{ij}(\cdot)$.  We now introduce more general functions that encompass all outlined approximations.  We begin with the general loss function $\check{L}_{\ell}(\cdot)$ that encompasses each functional form outlined in equations (\ref{LossFuncDef})-(\ref{AnotherApprox}). This general function highlights the convexity of the loss function around the origin as well as its symmetry about the point $\psi_{\ell}$.  Resistances are realistically assumed to be positive, resulting in strict convexity.  Note that a function with a check mark $ \ \check{ } \  $ represents an approximation to its \emph{exact} counterpart denoted with a hat symbol $ \ \hat{ } \  $.
\Def{\normalfont   A \emph{general loss function} of angles, denoted $\check{L}_{\ell}:\mathbb{R}\rightarrow\mathbb{R}$, is a function with the following properties: strictly convex on the subdomain $\mathcal{D}_{\ell}:=[-\frac{\pi}{2}+\psi_{\ell},\frac{\pi}{2}+\psi_{\ell}]$, continuously differentiable, non-negative, symmetric about the point $\psi_{\ell}$ and strictly monotonically increasing on the subdomain \mbox{$\mathcal{D}_{\ell +}:=[\psi_{\ell},\frac{\pi}{2}+\psi_{\ell}]$}.  \label{Lcheckdef}}\normalfont \vspace{-0pt}

The general mid-line power flow function $\check{F}_{\ell}(\cdot)$ encompasses all functional forms outlined by equations (\ref{midlinefloweq})-(\ref{midlinefloweqfinal}). This general function highlights the monotonic property of the mid-line power flow function for angle differences near the origin.  Although the mid-line power flow function is typically monotonically increasing this definition is left general to accommodate potentially positive susceptance values (or equivalently negative reactance values), in which case the function would be monotonically decreasing. 

\Def{\normalfont  A \emph{general mid-line power flow function} of angles, denoted $\check{F}_{\ell}:\mathbb{R}\rightarrow\mathbb{R}$, is strictly monotonic on the subdomain $\mathcal{D}_{\ell}$ and is continuously differentiable.\label{Fcheckdef}}\normalfont 

The general squared current magnitude function $\check{\mathcal{I}}_{ij}(\cdot)$ encompasses all functional forms outlined by equations (\ref{ExactCurrent})-(\ref{ISfinal}).  This general function highlights the convexity of the squared current magnitude function near the origin as well as its symmetry about the point $-\phi_{\ell}+\psi_{\ell}$.

\Def{\normalfont   A \emph{general squared current magnitude function} of angles, denoted $\check{\mathcal{I}}_{ij}:\mathbb{R}\rightarrow\mathbb{R}$, is a function with the following properties: convex on the subdomain \mbox{$\tilde{\mathcal{D}}_{\ell}:= [-\frac{\pi}{2}-\phi_{\ell}+\psi_{\ell},\frac{\pi}{2}-\phi_{\ell}+\psi_{\ell}]$}, strictly monotonically increasing on the subdomain \mbox{$\tilde{\mathcal{D}}_{\ell +}:= [-\phi_{\ell}+\psi_{\ell},\frac{\pi}{2}-\phi_{\ell}+\psi_{\ell}]$}, symmetric about the point $-\phi_{\ell}+\psi_{\ell}$ and continuously differentiable.  \label{Icheckdef}}\normalfont

The remainder of the paper provides results using these more general functions and is therefore pertinent to typical dispatch formulations that utilize the simplest functional forms outlined by equations (\ref{AnotherApprox}), (\ref{midlinefloweqfinal}), and  (\ref{ISfinal}) and also pertinent to the exact functional forms outlined by equations (\ref{LossFuncDef}), (\ref{midlinefloweq}), and (\ref{ExactCurrent}).  

\vspace{-5pt}
\section{Transmission Constrained Economic Dispatch}\label{sec:OPFForm}
\vspace{-5pt}
Economic Dispatch problems are solved to determine the optimal dispatched generation and set the locational marginal price.  In this section we formulate a non-convex economic dispatch problem using the general transmission line model from the previous section.  This problem formulation enforces line limits using simple bounds on the voltage angle difference across each transmission line.  We then explain how the voltage angle difference bounds can be chosen to enforce limits on a variety of line related quantities.

\vspace{-5pt}
\subsection{Power Injections}
\vspace{-5pt}
The net power injections at each bus can be expressed by summing the associated extractions from each transmission line adjacent to it.  To express the power injections as a function of voltage angles we will utilize the FND model (\ref{FNDEQ1}) and (\ref{FNDEQ2}) evaluated at the midpoint of each line $d=\frac{1}{2}$ along with arbitrary approximations for the loss function and mid-line power flows outlined in section \ref{sec:LineModel}.  Let's introduce a vector valued function $L: \mathbb{R}^m \rightarrow \mathbb{R}^m$ that maps voltage angle differences to line losses.  The $\ell^{th}$ element of this vector valued loss function is defined as $L_{\ell}(\Delta\!\theta):=\check{L}_{\ell}(\Delta\!\theta_{\ell})$ where $\Delta\!\theta\in \mathbb{R}^m$ is a vector of voltage angle differences, consistent with the established notation.  Similarly, a vector valued function $F: \mathbb{R}^m \rightarrow \mathbb{R}^m$ maps voltage angle differences to mid-line power flows.  The $\ell^{th}$ element of this vector valued function is defined as $F_{\ell}(\Delta\!\theta):=\check{F}_{\ell}(\Delta\!\theta_{\ell})$.

The net power injections at each bus in the system can be represented as a vector valued function of voltage angle differences, denoted $T: \mathbb{R}^{m}\rightarrow \mathbb{R}^{n}$:
\begin{equation} T(\Delta\!\theta):={\textstyle\frac{1}{2}}\left|A\right|^{\dag} L(\Delta\!\theta) + A^{\dag} F(\Delta\!\theta)  \label{PIdef}\text{,}\end{equation}
where the branch-bus incidence matrix of the graph $\mathcal{G}$ is denoted $A\in \mathbb{R}^{m \times n}$. Specifically, $A$ is sparse and the row representing line $\ell$ connecting bus $i$ to bus $j$ has element $i$ equal to $1$ and $j$ equal to $-1$. The element-wise absolute value of the branch-bus incidence matrix is denoted $\left|A\right|$, also known as the unoriented incidence matrix of graph $\mathcal{G}$.  Intuitively, the function $T(\cdot)$ can be interpreted as each bus being subject to a lossless power flow extraction associated with each line and also being assigned half the losses of each line incident to it. 

Since the voltage angles $\theta\in \mathbb{R}^n$ only enter each equation through the differences $\Delta\!\theta_{\ell}=A_{\ell}\theta$, one degree of freedom can be removed by assigning an arbitrary reference bus $\nu \in \mathbb{N}$ and setting $\theta_{\nu}=0$.  The vector of voltage angle differences can now be written as $\Delta\!\theta=\dot{A}\dot{\theta}$ where $\dot{\theta}\in \mathbb{R}^{n-1}$ is the vector of voltage angles with element $\nu$ removed and the matrix $\dot{A}\in \mathbb{R}^{m \times (n-1)}$ is the incidence matrix $A$ with column $\nu$ removed.

\subsubsection{Loss Distribution Factors}
Reference \cite{litvinov2004marginal} introduced loss distribution factors to distribute losses throughout the transmission system.  Loss distribution factors can be thought of as an approximation to the function of power injections $T(\Delta\!\theta)$.  A vector of distribution factors denoted $\eta\in \mathbb{R}^n_+$ sum to one and represent the fraction of total system losses allocated to each bus.  The associated approximation to the function representing power injections is written as follows:
\begin{equation} T(\Delta\!\theta)\approx \eta 1^{\dag} L(\Delta\!\theta) + A^{\dag} F(\Delta\!\theta)  \label{LDFs}\text{.}\end{equation}
Notice that $\eta 1^{\dag} \neq \frac{1}{2}\left|A\right|^{\dag}$ for any choice of loss distribution factors since the rank of $\eta 1^{\dag}$ is one and the rank of the unoriented incidence matrix $\left|A\right|$ is at least $n\hspace{-2pt}-\hspace{-2pt}1$ under the assumption that the system graph is fully connected~\cite{van1976incidence}.  Thus this approximation is never an exact representation of the function $T(\Delta\!\theta)$.  Instead, the proper choice of loss distribution factors will change with state $\dot{\theta}$.  In future work the authors will address the proper choice of loss distribution factors using linearization techniques.

\vspace{-5pt}
\subsection{Economic Dispatch Problem}\label{sec:ED}
\vspace{-5pt}
The economic dispatch problem optimizes over the nodal generation represented by vector $P\hspace{-3pt}\in\hspace{-3pt} \mathbb{R}^n$.  The cost of generation is represented by the function $C(P)$ where $C\hspace{-5pt}:\hspace{-5pt}\mathbb{R}^n\hspace{-5pt}\rightarrow \hspace{-5pt}\mathbb{R}$ is assumed convex.  The nodal demand is considered constant and is represented by $D\hspace{-2pt}\in\hspace{-2pt} \mathbb{R}^n$. The \emph{economic dispatch problem} is: 

\vspace{-10pt}
\small{
\begin{subequations}\label{OPF}
\begin{equation}\underset{P\in \mathbb{R}^n, \dot{\theta} \in \mathbb{R}^{n-1}}{\text{min}}  \ \ \ \ \   C(P) \tag{\ref{OPF}}\end{equation}
\vspace{-5pt}
\begin{equation}st:    T(\dot{A}\dot{\theta})=P-D  \label{OPFPowerBallance}\end{equation}
\vspace{-12pt}
\begin{equation}
 \ubar{P}\leq P\leq \bar{P} \label{OPF2}
\end{equation}
\vspace{-12pt}
\begin{equation}
 \underline{\Delta\!\theta} \leq \dot{A}\dot{\theta}  \leq \overline{\Delta\!\theta} \label{OPFPfeas}
\end{equation}

\end{subequations}}\normalsize

Constraint (\ref{OPFPowerBallance}) is a vector equality constraint that represents power balance at each node. Constraints (\ref{OPF2}) enforce generator output limits and constraints (\ref{OPFPfeas}) represent limits on voltage angle differences across each transmission line.  The next section will explain how to choose the voltage angle difference limits to indirectly enforce limits on various line related quantities including line losses and current magnitude.

Note that there is no explicit slack bus in this formulation, since power balance is represented at each bus, so the formulation is independent of choice of slack bus.  There is an explicit reference bus and Theorem \ref{ReferenceBusInvariance} below shows that the formulation is independent of choice of reference bus.

\begin{Thm}\label{ReferenceBusInvariance}Consider two instances of the economic dispatch problem (\ref{OPF}) defined using different reference buses.  The reduced voltage angle vectors associated with the first and second instances of the problem are defined to be $\dot{\theta}\hspace{-2pt}\in\hspace{-2pt} \mathbb{R}^{n-1}$ and $\ddot{\theta}\hspace{-2pt}\in\hspace{-2pt} \mathbb{R}^{n-1}$ respectively.  Let $(P^{\star}\hspace{-2pt},\dot{\theta}^{\star})$ be a solution to the first instance of the problem.  There exists some $\ddot{\theta}^{\star}\hspace{-4pt}\in\hspace{-2pt}\mathbb{R}^{n-1}$ such that $(P^{\star}\hspace{-2pt},\ddot{\theta}^{\star})$ is a solution to the second instance of the problem. 
\end{Thm}
\vspace{-5pt}
\proofname: 
Let $\dot{A}$ and $\ddot{A}$ be the reduced brach-bus incidence matrices for the first and second instances of the problem respectively. The matrix $\dot{A}$ and $A$ have the same range space because any given row of $A$ can be written as a linear combination of the other rows in $A$ (ie. $A1=0$). Similarly $\ddot{A}$ and $A$ have the same range space.  Thus $\dot{A}$ and $\ddot{A}$ have the same range space. As a result, there must exist some $\ddot{\theta}^{\star}$ such that $\ddot{A}\ddot{\theta}^{\star}=\dot{A}\dot{\theta}^{\star}$. It follows that there exists a $\ddot{\theta}^{\star}$ such that $(P^{\star},\ddot{\theta}^{\star})$ is feasible for the second instance of the problem. Notice that this implies the second instance of the problem has an optimal value no greater than the first instance.  Furthermore, there does not exist a feasible point of the second instance of the problem that achieves a lower optimal value than the first instance, else a feasible point can be constructed for the first instance that has a lower cost than the solution.  Thus both instances have the same optimal value and there exists a $\ddot{\theta}^{\star}$ such that $(P^{\star},\ddot{\theta}^{\star})$ is optimal for the second instance of the problem. \qed

The Locational Marginal Prices (LMPs), denoted $\lambda\in\mathbb{R}^n$, are defined to be the sensitivity of the optimal value of the economic dispatch problem (\ref{OPF}) with respect to the demand vector $D$, assuming such sensitivities exist.  In fact, the LMPs are independent of the choice of reference bus because the optimal value of the economic dispatch problem (\ref{OPF}) is independent of the choice of reference bus for any $D\in\mathbb{R}^n$ as stated in Theorem \ref{ReferenceBusInvariance}. Of course this simple analysis assumes that LMPs are well defined, which is a technicality that will be carified in future work.  Additional future work will decompose the LMP into energy, loss and congestion components.  Similar issues have been addressed by previous work \cite{moon2000slack}.

\remark{\normalfont A clear relationship must be established between the economic dispatch problem and the gold standard AC OPF problem, as specified in \cite{cain2012history}. Without explicitly defining the AC OPF problem, we note that it includes real and reactive power flows on transmission lines as decision variables and specifies fixed real and reactive power demand at buses.  In fact, under the assumption that shunt conductances are negligible, the \emph{exact} economic dispatch problem defined by equations (\ref{LossFuncDef}), (\ref{midlinefloweq}), and (\ref{ExactCurrent}) can be derived from the AC OPF problem by fixing the voltage magnitudes and removing all constraints that involve reactive power quantities, including reactive power balance constraints at buses.  Thus, if the voltage magnitudes are fixed to values that match the solution of the AC OPF problem, then the economic dispatch problem formulated here will act as a relaxation of the AC OPF problem.  Of course, in practice the solution to the AC OPF problem will not be available when fixing the voltage magnitudes, so these quantities must be approximated. \label{ACDecouple}}\normalfont
\vspace{-5pt}

\subsection{Characterizing Line Limits}\label{Sec:LineLim}

The economic dispatch problem should enforce line limits that represent the physical abilities of the transmission line.  For this reason it is appropriate to limit a real-power quantity or current quantity using approximations outlined in the previous section. Unfortunately, such constraints are non-linear and may even be non-convex.  This subsection explains how to enforce such limits by reformulating them to be in the standard form of constraints (\ref{OPFPfeas}), which place bounds on the voltage angle difference across each transmission line. 

We will consider two types of line limits. Namely, real-power flow limits enforced at the mid-point of the transmission line and limits on the squared current magnitude flowing into either side of the transmission line.  This framework can also accommodate other types of line limit constraints that are omitted due to space constraints.  Alternative line limit constraints include limits on the real-power loss across the transmission line and limits on the real-power injected into either end of the transmission line.  By reformulating these line limits to the same form as constraints (\ref{OPFPfeas}) this section implies that the economic dispatch problem (\ref{OPF}) encompasses all such limits without loss of generality.  


\remark{\normalfont Throughout the paper we will assume that any vector of voltage angle differences $\Delta\!\theta$ satisfying (\ref{OPFPfeas}) also satisfies the constraint $\Delta\!\theta_{\ell}\in \mathcal{D}_{\ell}\cap \tilde{\mathcal{D}}_{\ell}$ for each line $\ell$, where $\mathcal{D}_{\ell}$ and $\tilde{\mathcal{D}}_{\ell}$ are from definitions \ref{Lcheckdef} and \ref{Icheckdef}.  This assumption effectively enforces limits on the voltage angle differences $\Delta\!\theta$ that should hold for any power system. \label{SmallAngleAssumption}}\normalfont

\subsubsection{Mid-Line Power Flow Limits}\label{Sec:MidLim}
First consider power flow limits enforced at the mid-point of the transmission line.  Interpreted as lossless power flow limits in the FND model from section \ref{Sec:FND}, these constraints are written as follows.
\begin{equation}
	-\bar{F}_{\ell} \leq \check{F}_{\ell}(\Delta\!\theta_{\ell})\leq \bar{F}_{\ell}\label{midlims}
\end{equation}
We assume $\bar{F}_{\ell}$ is in the image of $\check{F}_{\ell}(\cdot)$ on the subdomain $\mathcal{D}_{\ell}$, denoted $\check{F}_{\ell}[\mathcal{D}_{\ell}]$.  As a result, limits on the mid-line power flow are easy to enforce because the function $\check{F}_{\ell}(\cdot)$ is strictly monotonic on the specified subdomain.  Define the function $\check{F}_{\ell}^{-1}:\check{F}_{\ell}[\mathcal{D}_{\ell}]\rightarrow \mathcal{D}_{\ell}$ as the inverse of $\check{F}_{\ell}(\cdot)$ on the specified subdomain.  This inverse function is strictly monotonic on $\check{F}_{\ell}[\mathcal{D}_{\ell}]$, allowing constraints (\ref{midlims}) to be written as follows.  Notice that these constraints are in the same form as (\ref{OPFPfeas}).
\begin{equation}
	\check{F}^{-1}_{\ell}(-\bar{F}_{\ell}) \leq \Delta\!\theta_{\ell}\leq \check{F}^{-1}_{\ell}(\bar{F}_{\ell})\label{midlims2}
\end{equation}

\subsubsection{Current Magnitude Limits}\label{CurrMagLimSec}
It may be more appropriate to limit the magnitude of the current flowing into the transmission line.  A limit on the squared magnitude of the current flowing into the line is represented by the following constraint.
\begin{equation}\check{\mathcal{I}}_{ij}(\Delta\!\theta_{\ell})\leq \bar{\mathcal{I}}^2_{\ell} \label{currconst}\end{equation}
where $\bar{\mathcal{I}}_{\ell}$ is the constant current magnitude limit and $\bar{\mathcal{I}}^2_{\ell}$ is assumed to lie in the image of $\check{\mathcal{I}}_{ij}(\cdot)$ on the subset $\tilde{\mathcal{D}}_{\ell}$, denoted $\check{\mathcal{I}}_{ij}[\tilde{\mathcal{D}}_{\ell}]$. Notice that $\check{\mathcal{I}}_{ij}(\cdot)$ is symmetric about the point $-\phi_{\ell}+\psi_{\ell}$ and is strictly monotonically increasing on the domain $\tilde{\mathcal{D}}_{\ell +}$ and so it is invertible on this subdomain.  Define $\check{\mathcal{I}}_{ij}^{-1}:\check{\mathcal{I}}_{ij}[\tilde{\mathcal{D}}_{\ell +}]\rightarrow \tilde{\mathcal{D}}_{\ell +}$ to be the inverse of the function $\check{\mathcal{I}}_{ij}(\cdot)$ restricted to the subdomain $\tilde{\mathcal{D}}_{\ell +}$. This inverse function is also strictly monotonically increasing on $\check{\mathcal{I}}_{ij}[\tilde{\mathcal{D}}_{\ell +}]$, allowing constraint (\ref{currconst}) to be rewritten as constraint (\ref{lossconst2}). This constraint can be easily placed in the same form as (\ref{OPFPfeas}). 
\begin{equation}| \Delta\!\theta_{\ell}+\phi_{\ell}-\psi_{\ell}|-\phi_{\ell}+\psi_{\ell}\leq \check{\mathcal{I}}^{-1}_{ij}(\bar{\mathcal{I}}^2_{\ell}) \label{lossconst2}\end{equation}

\remark{\normalfont Since voltage magnitudes are assumed fixed, the magnitude of current flowing into either side of a transmission line is proportional to the magnitude of the apparent power flowing into that side of the transmission line.  Thus limits on the magnitude of apparent power flow (MVA limits) can be enforced using squared current magnitude limits (\ref{currconst}) and (\ref{lossconst2}).}\normalfont


\vspace{-10pt}
\section{Empirical Results} \label{sec:Results}
\vspace{-5pt}

This section provides an empirical analysis of three realistically large test cases provided by version 6.0 of the MATPOWER toolbox in MATLAB \cite{zimmerman2011matpower}.  The \emph{3375wp} test case is a 3,375 bus representation of the Polish power system during the winter 2007-2008 winter evening peak with a total fixed demand of $48,362$ MW. The \emph{2869pegase} test case is a 2,869 bus representation of the European high voltage transmission network with a total fixed demand of $132,437$ MW. The \emph{6515rte} test case is a 6,515 bus representation of the French transmission network with a total fixed demand of $107,264$ MW.  The two larger test cases are fully described in reference \cite{josz2016ac}. For each test case the MVA rating of each transmission line is interpreted as a current magnitude limit and line limits are enforced as in section \ref{CurrMagLimSec}. 


%

We consider six different optimization problems for each test case: the AC OPF problem, the exact economic dispatch problem and four approximations to the exact economic dispatch problem. Each optimization problem is solved by the interior-point algorithm provided by the MATLAB function FMINCON with user supplied analytical gradients/hessians using a standard laptop with a 2.7 GHz processor. Due to the non-convexity of each problem, global optimality cannot be guaranteed in general and the interior point algorithm may converge to a local minimizer.  That being said, we are able to verify that a global minimizer was identified for each of the four approximations for test cases \emph{6515rte} and \emph{2869pegase} by use of the load over-satisfaction relaxation described later in this section.

The AC OPF problem fully captures the coupling of real and reactive power and optimizes over voltage magnitudes. At the optimal dispatch of the AC OPF problem, the operating cost of test cases \emph{6515rte}, \emph{2869pegase}, and \emph{3375wp} are \$$109,767$, \$$133,993$, and \$$7,404,635$ respectively.  The \emph{exact} economic dispatch problem (\ref{OPF}) is defined by the exact expressions of each function provided in the body of the paper (\ref{LossFuncDef}), (\ref{midlinefloweq}), (\ref{ExactCurrent}) and (\ref{PIdef}) and is formulated by fixing voltage magnitudes to the identified local minimizer of the AC OPF problem. As described in remark \ref{ACDecouple}, this problem should act as a relaxation and an approximation to the AC OPF problem. The identified local minimizer of the AC OPF problem is used to initialize the interior point algorithm for the exact economic dispatch problem.  All values associated with the identified local minimum of this problem are denoted with a hat $\hat{ }$ \hspace{1pt}.

\begingroup
\setlength{\tabcolsep}{1pt} 
\renewcommand{\arraystretch}{1} 

\begin{table*}[b!]

\vspace{-10pt}
\caption{\label{LargeSysResults} This table compares the solution of multiple different economic dispatch problems to that of the exact economic dispatch problem using three realistically large test cases.  Quantities associated with the exact formulation are denoted by a hat and quantities associated with general formulations are denoted by a $ \ \check{ } \  $. The reported time is an average over $10$ runs.}
\vspace{-10pt}
\label{tt}
\begin{center}
\begin{tabular}{|c||c|c||c||c|c|c|c||c|c||c|c|c|} \hline
Test &Aprx.&Aprx. Eqns.&$C(\hat{P})\hspace{-3pt}-\hspace{-3pt}C(\check{P})$ & $1^T \check{P}$ &$\|\check{P}\|_0$& $\|\check{P}\hspace{-3pt}-\hspace{-3pt}\hat{P}\|_{1}$&$\|\check{P}\hspace{-3pt}-\hspace{-3pt}\hat{P}\|_{\infty}$& \rule{0pt}{9pt} mean($\check{\lambda}$)&min($\check{\lambda}$) / max($\check{\lambda}$) &\multirow{ 2}{*}{$\underset{i\in\mathbb{N}}{\text{max}} \ \frac{|\hat{\lambda}_i\hspace{-2pt}-\hspace{-2pt}\check{\lambda}_i|}{\rule{0pt}{7pt}\hat{\lambda}_i}$}&$\|\check{\Delta\!\theta}\|_{\infty}$&time\\
Case &Num.&$\check{F}_{\ell}$ / $\check{L}_{\ell}$ / $\check{\mathcal{I}}_{ij}$ / $T$& (\$) &  (MW) & (MW) & (MW) & (MW)&(\$/MW) &(\$/MW)&&(rad)&(sec)\\ \hline \hline  
     
\multirow{ 6}{*}{\hspace{-5pt}6515rte\hspace{-5pt}}& Exact &(\ref{LossFuncDef}) / (\ref{midlinefloweq}) / (\ref{ExactCurrent}) / (\ref{PIdef})& 0 & 109764 & 652 & 0 & 0 & 1.08 & 0.87 / 2.52 & 0 & 0.3822 & 10\\ \cline{2-13}
     & 1 &(\ref{LossFuncFirstApprox}) / (\ref{midlineflowfirstapprox}) / (\ref{ISFirstApprox}) / (\ref{PIdef})&2.36 & 109762 & 643 & 159.29 & 31.06 & 1.08 & 0.87 / 2.53 & 0.0891 & 0.3730 & 8\\ \cline{2-13}
     & 2 &(\ref{approx1}) / (\ref{midlinefloweqtaylor}) / (\ref{IStaylor}) / (\ref{PIdef})&-184.57 & 109949 & 673 & 8622.52 & 658.90 & 1.07 & 0.97 / 1.43 & 0.4306 & 0.4255 & 8\\ \cline{2-13}
     & 3 &(\ref{AnotherApprox}) / (\ref{midlinefloweqfinal}) / (\ref{ISfinal}) / (\ref{PIdef})&-72.21 & 109837 & 677 & 8793.72 & 664.70 & 1.07 & 0.97 / 1.41 & 0.4394 & 0.4254 & 24\\ \cline{2-13}
     & 4 &(\ref{AnotherApprox}) / (\ref{midlinefloweqfinal}) / (\ref{ISfinal}) / (\ref{LDFs})&47.92 & 109717 & 677 & 8012.29 & 679.46 & 1.07 & 0.96 / 1.55 & 0.3829 & 0.4255 & 27\\\cline{2-13}
     & AC OPF & &-3.38 & 109768 & 643 & 908.16 & 103.87 & 1.08 & 0.56 / 2.48 & 0.5257 & 0.3822 & 152\\ \hline \hline 
     

  \multirow{ 6}{*}{2869pegase} & Exact &(\ref{LossFuncDef}) / (\ref{midlinefloweq}) / (\ref{ExactCurrent}) / (\ref{PIdef})& 0 & 133982 & 468 & 0 & 0 & 1.02 & 0.99 / 1.12 & 0 & 0.2429 & 4\\ \cline{2-13}
     & 1 &(\ref{LossFuncFirstApprox}) / (\ref{midlineflowfirstapprox}) / (\ref{ISFirstApprox}) / (\ref{PIdef})&2.14 & 133980 & 468 & 58.06 & 5.21 & 1.02 & 0.99 / 1.12 & 0.0011 & 0.2422 & 3\\ \cline{2-13}
     & 2 &(\ref{approx1}) / (\ref{midlinefloweqtaylor}) / (\ref{IStaylor}) / (\ref{PIdef})&-167.57 & 134150 & 487 & 2124.47 & 190.99 & 1.02 & 0.98 / 1.14 & 0.0330 & 0.2780 & 3\\ \cline{2-13}
     & 3 &(\ref{AnotherApprox}) / (\ref{midlinefloweqfinal}) / (\ref{ISfinal}) / (\ref{PIdef})&-134.32 & 134117 & 487 & 1971.91 & 188.48 & 1.02 & 0.98 / 1.13 & 0.0177 & 0.2743 & 7\\ \cline{2-13}
     & 4 &(\ref{AnotherApprox}) / (\ref{midlinefloweqfinal}) / (\ref{ISfinal}) / (\ref{LDFs})&-94.62 & 134077 & 494 & 3739.75 & 972.89 & 1.02 & 0.98 / 1.12 & 0.0132 & 0.2710 & 10\\\cline{2-13}
     & AC OPF & &-11.29 & 133993 & 460 & 581.57 & 22.25 & 1.02 & 0.99 / 1.11 & 0.0222 & 0.2438 & 6\\ \hline \hline

     \multirow{ 6}{*}{3375wp} & Exact &(\ref{LossFuncDef}) / (\ref{midlinefloweq}) / (\ref{ExactCurrent}) / (\ref{PIdef})& 0 & 49190 & 469 & 0 & 0 & 146.70 & 0.00 / 454.89 & 0 & 0.2639 & 8\\ \cline{2-13}
     & 1 &(\ref{LossFuncFirstApprox}) / (\ref{midlineflowfirstapprox}) / (\ref{ISFirstApprox}) / (\ref{PIdef})&344.10 & 49188 & 469 & 32.79 & 16.58 & 146.66 & 0.00 / 454.58 & 0.9604 & 0.2624 & 5\\ \cline{2-13}
     & 2 &(\ref{approx1}) / (\ref{midlinefloweqtaylor}) / (\ref{IStaylor}) / (\ref{PIdef})&-33629.36 & 49266 & 469 & 2267.75 & 309.97 & 151.26 & 0.00 / 971.06 & 5.7514 & 0.2625 & 6\\ \cline{2-13}
     & 3 &(\ref{AnotherApprox}) / (\ref{midlinefloweqfinal}) / (\ref{ISfinal}) / (\ref{PIdef})&-13737.40 & 49228 & 469 & 1011.08 & 245.73 & 149.27 & 0.00 / 460.28 & 1.6180 & 0.2726 & 6\\ \cline{2-13}
     & 4 &(\ref{AnotherApprox}) / (\ref{midlinefloweqfinal}) / (\ref{ISfinal}) / (\ref{LDFs})&-7382.55 & 49184 & 469 & 911.03 & 244.28 & 148.55 & 0.00 / 459.67 & 539.2388 & 0.2612 & 7\\\cline{2-13}
     & AC OPF & &-8.43 & 49188 & 469 & 118.84 & 54.03 & 147.42 & 0.00 / 469.40 & 0.9943 & 0.2629 & 180\\
     \hline 
\end{tabular}
\end{center}
\vspace{-17pt}
\end{table*}

\endgroup

A relaxed version of the general economic dispatch problem that uses the \emph{load-over satisfaction relaxation} can be obtained from problem (\ref{OPF}) by replacing the power balance equality constraints (\ref{OPFPowerBallance}) with inequality constraints that permit the delivery of excess demand~\cite{palma2013modelling}.  If this \emph{relaxed economic dispatch problem} results in positive LMPs, then its solution also solves the general economic dispatch problem (\ref{OPF}).  Additionally, the relaxed economic dispatch problem can be proven convex under the condition that an affine approximation of the mid-line power flow function is used and thus the proposed interior point methods are guaranteed to converge to a global minimizer.  Fortunately, all proposed approximations to the mid-line power flow function in section \ref{MidLinePowerFlowAproximations} satisfy this property.  This allows us to verify that the interior point methods identified global optima for test cases \emph{6515rte} and \emph{2865pegase} when using each of the four approximations to the exact economic dispatch problem.  However, we emphasize that this approach cannot always be taken in the case where there exists an LMP that is non-positive as in test case \emph{3375wp}, see reference~\cite{philpott2004financial}. We found no large test cases with any strictly negative LMPs.

The remainder of this section analyzes table \ref{LargeSysResults}, which provides detailed information about the six aforementioned optimization problems solved for each test case.  The identified local minimizers of the AC OPF problem and the approximations 1-4 are directly compared to that of the exact economic dispatch problem. All values pertaining to the local minimizers of a general optimization problem are denoted with a check mark, ie. the dispatch is $\check{P}$.  Approximation 1 uses Taylor expansions to obtain a very accurate quadratically constrained program.  Approximation 2 additionally assumes voltage magnitudes are nominal, tap ratios are nominal and shunt susceptances are much smaller than series susceptances.  Approximation 3 additionally assumes that series resistances are much smaller than series reactances.  Approximation 4 additionally uses the load distribution factor approximation with the LDFs chosen to allocate all losses to the slack bus, which is designated by each individual test case description.  The third column of the table explicitly states the equations used in each approximation.

\vspace{-3pt}
\remark{\normalfont From table \ref{LargeSysResults}, the voltage angle difference across each line falls well within the limits from remark \ref{SmallAngleAssumption}, which approximately constrains $|\Delta\!\theta_{\ell}|$ to be lower than $\frac{\pi}{2}\hspace{-2pt}\approx\hspace{-2pt} 1.5707$.} \normalfont


\vspace{-10pt}
\subsection{Economic Dispatch Problem vs. AC OPF}
\vspace{-5pt}
We begin by comparing the economic dispatch problem (\ref{OPF}) to the AC OPF problem provided by MATPOWER.  Each test case follows the trend outlined in remark \ref{ACDecouple}.  Specifically, the \emph{exact} economic dispatch problem acts as a relaxation to the AC OPF problem as it attains a lower optimal value.  However, the optimal values of both problems are relatively close to each other as compared to the total system cost. 


It is perhaps more important to analyze the difference in the generation dispatch between the exact economic dispatch problem and the AC OPF problem.  The 1-norm of the nodal dispatch approximation error, denoted $\|\check{P}\hspace{-2pt}-\hspace{-2pt}\hat{P}\|_1$, is relatively small, realizing values of no more than $1\%$ of total system demand.  In contrast, the infinity norm of the nodal dispatch error is potentially significant because it represents a MW value that is seen entirely by a single bus.  In fact, marginal generators see the largest change in dispatch.  The sum of all dispatched generation is denoted $1^T\check{P}$ and is similar across all approximations.  The zero-norm of the nodal dispatch, denoted $\|\check{P}\|_0$, represents the number of non-zero elements in the vector $\check{P}$.  This quantity represents the number of dispatched generators and does not change much between approximations. 

A similar conclusion can be drawn for the LMPs identified by both problems.  In general LMPs are very similar for both problems, illustrated by nearly identical mean LMPs.  However, the LMP at an individual bus may be significantly different for both problems as illustrated by the maximum normalized difference of nodal price denoted $\frac{|\hat{\lambda}_i\hspace{-2pt}-\hspace{-2pt}\check{\lambda}_i|}{\check{\lambda_i}}$.  Although this value is very small for test case \emph{2869pegase}, at least one bus in test case \emph{3375wp} has an LMP change of over 99\% in magnitude.

\vspace{-5pt}
\subsection{Approximations of the Exact Economic Dispatch Problem}
\vspace{-5pt}

Table \ref{LargeSysResults} quantifies how well each economic dispatch formulation approximates the exact formulation.  Approximation 1 is the most accurate in terms of nodal dispatch and LMP approximation error.  This is expected because only accurate Taylor expansion approximations are used.  

Approximation 2 experiences a drastic increase in approximation error because it introduces multiple assumptions including nominal tap ratios, negligible shunt susceptances and nominal voltage magnitudes.  

Approximation 3 additionally assumes $r_{\ell} \hspace{-3pt} \ll \hspace{-3pt} x_{\ell}$ and the resulting approximate loss function tends to underestimate real-power losses as compared to approximation 2.  This can be seen by noticing that the total dispatched generation is lower for approximation 3.  

Approximation 4 is equivalent to approximation 3 but with all losses allocated to the slack bus. Notice that the nodal dispatch error has a large infinity norm, $\|\check{P}\hspace{-2pt}-\hspace{-2pt}\hat{P}\|_{\infty}$, because the load profile changes significantly.  For this reason approximation 4 can be very inaccurate.  For example the 1-norm of the nodal dispatch error is nearly doubled in test case \emph{2869pegase}.  

\vspace{-5pt}
\section{Conclusions} \label{sec:Conc}
\vspace{-5pt}
This paper derives a generalized non-convex economic dispatch problem with marginal losses that consistently models flow and loss approximation, results in approximately correct outcomes and is proven to be reference bus independent.  A hierarchy of approximations for this problem are outlined and common unnecessary assumptions are identified along with proper choices of ``tuning parameters.''  For example, nodal loss allocation is derived from first principles to assign half losses of each line to its adjacent buses.  Additionally, line limit constraints are derived in the form of current magnitude limits and mid-line power flow limits and are enforced as simple bounds on the voltage angle difference on each transmission line.  Empirical results are provided that illustrate the general trend of increasing approximation error as more approximations are used.  Furthermore, the identified local minimizer of the exact economic dispatch problem is shown to very closely match the certified global minimizer of an approximate economic dispatch problem solved via the load over-satisfaction relaxation. Certain approximations increase the error in price and dispatch outcomes; however, the approximation that allocates losses to the slack bus is very inaccurate and results in significant dispatch and LMP errors.


\vspace{-10pt}

\bibliography{bibfile}
\bibliographystyle{IEEEtran}

\end{document}